\documentclass[conference]{IEEEtran}
\IEEEoverridecommandlockouts
\usepackage{cite}
\usepackage{amsmath,amssymb,amsfonts}
\usepackage{algpseudocode}
\usepackage{algorithm}
\usepackage{algorithmicx}
\usepackage{graphicx}
\usepackage{textcomp}
\usepackage{xcolor}
\usepackage{hyperref}
\usepackage{soul}
\def\BibTeX{{\rm B\kern-.05em{\sc i\kern-.025em b}\kern-.08em
    T\kern-.1667em\lower.7ex\hbox{E}\kern-.125emX}}

\newcommand{\eps}{\varepsilon}
\begin{document}

\title{Low synchronization GMRES algorithms\\
\thanks{This work was funded by the DOE Exascale Computing Project (17-SC-20-SC).}
}

\author{\IEEEauthorblockN{Katarzyna Swirydowicz, Julien Langou, Shreyas Ananthan, Ulrike Yang, Stephen Thomas}
\IEEEauthorblockA{ \footnotesize
National Renewable Energy Laboratory, Golden, Colorado, U.S.A.  \\
University of Colorado Denver, Colorado, U.S.A. \\
Lawrence Livermore National Laboratory, Livermore, California, U.S.A. \\
Katarzyna.Swirydowicz@nrel.gov, julien.langou@ucdenver.edu, Shreyas.Ananthan@nrel.gov, 
yang11@llnl.gov, Stephen.Thomas@nrel.gov
} }

\maketitle

\begin{abstract}
Communication-avoiding and pipelined variants of Krylov solvers are critical
for the scalability of linear system solvers on future exascale architectures.
We present low synchronization variants of iterated classical (CGS) and
modified Gram-Schmidt (MGS) algorithms that require one and two global
reduction communication steps.  Derivations of low synchronization iterated CGS
algorithms are based on previous work by Ruhe.  Our main contribution is to
introduce a backward normalization lag into the compact $WY$ form of MGS
resulting in a ${\cal O}(\eps)\kappa(A)$ stable GMRES algorithm that requires
only one global synchronization per iteration.  The reduction operations are
overlapped with computations and pipelined to optimize performance. Further
improvements in performance are achieved by accelerating GMRES BLAS-2
operations on GPUs.
\end{abstract}


\section{Introduction}

The Gram-Schmidt algorithm constructs an orthogonal basis from a set of
linearly independent vectors spanning a vector space. The classical (CGS)
algorithm is column-oriented, whereas the modified (MGS) Gram-Schmidt algorithm
has both row and column variants. The numerical properties of both variants were
clarified by Bj\"orck \cite{Bjorck67}, who established the backwards stability
of the $QR$ factorization of a matrix $A$ produced by MGS for least squares
problems. The columns of $Q$ lose othogonality in both cases and bounds on
$\|I - Q^TQ\|_2$ have been derived and improved by several authors.  For the
modified algorithm, Bj\"orck \cite{Bjorck67} derived a bound based on the
condition number $\kappa(A)$.  The situation is worse for the classical
algorithm where the loss of orthogonality depends on $\kappa^2(A)$ as shown in
\cite{Giraud05} and \cite{Giraud06}. For both CGS and MGS, the orthogonality
can be restored to the level of machine round-off by reorthogonalization of the
columns of $Q$. Iterative reorthogonalization was discussed in Leon
et.~al.~\cite{Leon10} along with the ``twice is enough" result of Kahan %\cite{Parlett98}, 
which was later generalized by Giraud et.~al.~\cite{Giraud06}.
A subsequent paper by Smoktunowicz et.~al.~ \cite{Smok06}
introduced a more stable Cholesky-like CGS algorithm using a Pythagorean
identity for the normalization of the basis vectors.

We present low synchronization variants of iterated classical (CGS) and
modified Gram-Schmidt (MGS) algorithms that require one and two global
reduction communication steps. The MGS algorithms are based on the compact $WY$
form of the orthogonalization process described in Bjorck \cite{Bjorck94} and
Malard and Paige \cite{Malard94}.  Our main contribution is to introduce a lag
into the column oriented MGS algorithm in order to delay normalization of the
diagonal elements of $R$.  In addition, the reductions can be overlapped with
the computation and pipelined as demonstrated by Ghysels et.~al.~\cite{Ghysels13} and
Yamazaki et.~al.~\cite{Yamazaki06}.  Derivations of low communication
algorithms for iterated CGS can be based on the paper by Ruhe \cite{Ruhe83} who
establishes the equivalence of reorthogonalization for classical Gram-Schmidt
with the Gauss-Jacobi iteration for the normal equations with coefficient
matrix $Q$ and Gauss-Seidel for modified Gram-Schmidt.  The author unrolls the
recurrence for re-orthogonalization, thus obtaining a loop-invariant form of
the algorithm.  The loop invariant allows us to combine the two
reorthogonalization iterations into one step and thereby eliminate the need for
an additional global reduction. The loss of orthogonality relationship for the
matrix $Q^TQ$ is also exploited to reduce the computation and communication overhead.

Both the classical and modified Gram-Schmidt algorithms may employ BLAS-2
vector inner products and norms to construct an orthonormal set of basis
vectors.  On modern multicore computers global reductions use a tree-like sum
of scalar products involving every core (MPI rank) on every node of the
machine.  The latency of this operation grows with the height of the tree
determined by the number of cores.  The remaining operations required in the
orthogonalization are scalable.  These consist of the mass inner-product $v =
X^Ty$, vector MAXPY operation $w = z - VH_{:,i}$, and sparse matrix-vector
product (SpMV) at each iteration.  For sparse matrices distributed by rows, the
parallel SpMV typically requires local communication between cores. The BLAS-1
vector sum is trivially data-parallel.  Hence, the global reductions are the
limiting factor for the available parallelism in the algorithm. Our MGS
algorithm reduces the number of global reductions per iteration to only one.

\section{Gram-Schmidt and Reorthogonalization}

Ruhe \cite{Ruhe83} provides a unique perspective by proving that the iterated
classical and modified Gram-Schmidt algorithms are equivalent to Gauss-Jacobi
and Gauss-Seidel iterations for linear systems. Let $A$ be an $n\times n$
non-singular matrix with $QR$ factorization obtained from the Gram-Schmidt
algorithm.  Consider one step of the classical Gram-Schmidt algorithm, with $p$
orthonormal vectors forming an $n\times p$ matrix $Q$, and the vector $a$,
which is to be orthogonalized to form the $p + 1$--st column of $Q$. The
classical Gram-Schmidt algorithm with re-orthogonalization is presented as
Algorithm \ref{alg:cgs}.
\begin{algorithm}
\footnotesize
\caption{\footnotesize Classical Gram-Schmidt Algorithm with re-orthogonalization}

\begin{algorithmic}[1]
  \State{Start $a^0 = a$, $r^0 = 0$.}
  \For {$k=1, 2, \ldots$}
     \State{$s^{k-1} = Q^Ta^{k-1}$ }
     \State{$r^k = r^{k-1} + s^{k-1}$ }
     \State{$a^k = a^{k-1} - Qs^{k-1}$ }
  \EndFor
  \State{$q_{p+1} = a^k / r_{p+1,p+1} = \|a^k\|_2$}
\end{algorithmic}
\label{alg:cgs}
\end{algorithm}
Each iteration in Steps 3--5 is a reorthogonalization. The notation in
Step 7 means that the final vector $a^k$ is normalized to form the next column
$q_{p+1}$ of $Q$. The normalization coefficient enters the $R$ matrix below
the vector $r$, leading to the corresponding matrix factors
\[
R_{p+1} =
\left[
\begin{array}{cc}
R_p & r \\
    & r_{p+1,p+1}
\end{array}
\right], \quad
Q_{p+1} =
\left[
\begin{array}{cc}
Q & q_{p+1}
\end{array}
\right]
\]
For the column-wise modified Gram-Schmidt algorithm, each iteration
(Steps 2 to 4) in Algorithm \ref{fig:MGSA} is divided into $p$ rank-1 updates, 
\begin{algorithm}
\footnotesize
\centering
\begin{algorithmic}[1]
  \For {$i=1, \ldots, p$}
     \State{$s_i^{k-1} = q_i^Ta^{k-1,i-1}$ }
     \State{$a^{k-1,i} = a^{k-1,i-1} - q_i s_i^{k-1}$ }
     \State{$r_i^k = r_i^{k-1} + s_i^{k-1}$ }
  \EndFor
  \State{$a_{k} = a^{k-1,p}$ }
\end{algorithmic}
\caption{\footnotesize \label{fig:MGSA} Modified Gram-Schmidt Algorithm }
\end{algorithm}

The key observation made by Ruhe \cite{Ruhe83} is that for the classical Gram-Schmidt algorithm, the
loop-invariants are obtained by unrolling the recurrences as follows,
\begin{eqnarray*}
a^k & = & a^{k-1} - Q s^{k-1} \\
    & = & a^{k-2} - Q s^{k-2} - Q s^{k-1} \\
    & = & a - Q (\: s^0 + s^1 + \cdots + s^{k-1} \:)
\end{eqnarray*}
\begin{eqnarray*}
r^k & = & r^{k-1} + s^{k-1} \\
    & = & r^{k-2} + s^{k-2} + s^{k-1} \\
    & = & r^0 + s^0 + s^1 + \cdots + s^{k-1}
\end{eqnarray*}
and therefore $a^k = a - Qr^k = Pa^{k-1}$, where $P$ is a projection matrix. It follows that
\begin{eqnarray*}
r^k & = & r^{k-1} + Q^T a^{k-1} \\
    & = & r^{k-1} + Q^T ( \: a - Qr^{k-1} \: )\\
    & = & r^{k-1} + Q^T  a - Q^TQ r^{k-1}
\end{eqnarray*}
This is the Gauss-Jacobi iteration applied to the normal equations
for the matrix $Q$ and right hand side vector $a$
\[
Q^TQ\:r = Q^T\:a
\]
The loop-invariant form of the matrix update can be judiciously applied in
order to minimize the number of communication steps in parallel
implementations of classical Gram-Schmidt.

For distributed-memory parallel computation based on message-passing, the
number of global reduction summation steps in the iterated classical
Gram-Schmidt algorithm can be reduced to two by employing the Ruhe
loop-invariant.  In particular, it is possible to write Step 4 of Algorithm~\ref{alg:cgs}
by combining the first two iterations of the re-orthogonalization as follows,
given $r^1 = Q^Ta$
\begin{eqnarray*}
r^2 & = & r^1 + Q^T a - Q^TQ r^1 \\
    & = & Q^Ta + Q^Ta - Q^TQ r^1 \\
    & = & Q^Ta + (\:I - Q^TQ\:) (\: Q^Ta \:)
\end{eqnarray*}
The above derivation implies that the projection matrix takes the form,
\[
P = I - Q(\: 2I - Q^TQ \:)Q^T = I - QTQ^T, \quad T = I - L - L^T
\]
where $L$ is a strictly lower triangular matrix
\[
Q^TQ = I + L + L^T
\]
The upper triangular matrix $L^T$ can be computed one column at a time and this suggests 
the form of a two synchronization step parallel algorithm for iterated CGS 
(Algorithm~\ref{fig:CGS2}).
\begin{algorithm}
\footnotesize
\centering
\begin{algorithmic}[1]
     \State{$[\:y,\: L_{:,p}^T\:]  = Q^T[\:a,\: q_p\:]$}
     \State{$r =  Ty$}
     \State{$a = a - Q Q^T a$}
     \State{$r_{p+1,p+1} = \|a\|_2$ }
     \State{$q_{p+1} = a/r_{p+1,p+1}$ }
\end{algorithmic}
\caption{\footnotesize \label{fig:CGS2} Iterated Classical Gram-Schmidt Algorithm with Two Synchronizations}
\end{algorithm}
Two global reduction communications are required for the mass inner products in Steps 1,
and the normalization in Step 4. The above algorithm suggests that a synchronization step
can be eliminated if the normalization Step 4 was lagged to the $p$-th iteration.
A lagged modified Gram-Schmidt algorithm based on the compact $WY$ representation
of the projection matrix $P$ is derived below.

In the case of the modified Gram-Schmidt algorithm, the loop invariant form
of the $(k-1)$-st step of the re-orthogonalization update for a column $i$
of the matrix $A$ is derived by Ruhe as

\[
a^{k-1, i} = a - Qr^{k-1, i}
\]
where 
\[
r^{k-1, i} = 
\left[ \begin{array}{cccccc}
r_1^k, & \ldots, & r_{i-1}^k, & r_i^{k-1}, & \ldots, & r_p^{k-1}
\end{array} \right]^T
\]
The columns of the matrix $R$ are updated as follows,
\[
r^k = (\:I + L\:)^{-1}Q^Ta - (\:I + L\:)^{-1} L^T r^{k-1},
\]

An important matrix power series expansion links the different Gram-Schmidt algorithms
to the stability analyses and the loss of orthogonality relations. In addition,
the equivalence of the modified Gram-Schmidt algorithm and compact $WY$ representation 
given below follow from the power series
\[
T = (\: I + L \:)^{-1} - I = -L + L^2 - L^3 + L^4 - \cdots
\]
A measure of the loss of orthogonality was introduced by Paige \cite{Paige18}
as $\|S\|_2$, where $S = (\:I + L^T\:)^{-1} L^T$. The norm remains close to
${\cal O}(\eps)$ for orthogonal vectors and increases to one as orthogonality
is lost.  However, given the loop-invariants, the matrix $L$ derived by Ruhe
does not contain the higher-order inner-products that are present in the $L_1$
matrix obtained for the compact $WY$ MGS projection.

In order to derive a block form of the modified Gram-Schmidt algorithm 
based on level-2 BLAS type operations, Bj\"orck \cite{Bjorck94} proposed
a compact $WY$ representation for MGS analogous to the Householder
factorization \cite{Scheiber89}. The matrix projection $P$ can be written as
\begin{eqnarray*}
P & = & (\: I - q_nq_n^T \:)\:(\: I - QL_1Q^T\:) \\
  & = & I - q_nq_n^T - QL_1Q^T + q_nq_n^TQL_1Q^T
\end{eqnarray*}
where $L_1$ is the lower triangular matrix of basis vector inner products. 
The $WY$ representation is then given by the matrix form
\[
P = I - 
\left[ \begin{array}{cc} Q & q_n \end{array} \right]
\left[ 
\begin{array}{cc} 
L_1       & 0 \\
-q_n^TQL_1 & 1 
\end{array} 
\right]
\left[ \begin{array}{c} Q^T \\ q_n^T \end{array} \right]
\]
and where
\[
P = I - QTQ^T, 
\]
The paper by Malard and Paige \cite{Malard94} shows how the transpose of the
matrix $T$ above can be formed recursively as follows
\begin{eqnarray*}
T_0 & = & 1\\
T_k & = &
\left[ 
\begin{array}{cc} 
T_{k-1} & -T_{k-1}Q^Tq_n \\
0       & 1 
\end{array} 
\right]
\end{eqnarray*}

A parallel modified Gram-Schmidt algorithm that requires only one global synchronization
step can be based on the above compact $WY$ form combined with a lagged normalization
of the diagonal of the upper triangular matrix $R$. 
\begin{algorithm}
\footnotesize
%\scriptsize
\centering
\begin{algorithmic}[1]
%\State{$\left[ \begin{array}{ccc} Q_{1:m,1:j}, & R_{1:j,1:j}, & T_{1:j,1:j} \end{array} \right] 
\State{${\rm mgs\_lvl2} 
\left[ \begin{array}{ccc} Q_{1:m,1:j}, &  R_{1:j,1:j}, & T_{1:j,1:j} \end{array} \right]$}
     \State{$ [\: Q_{j-1}^Tq_{j-1},\: Q_{j-1}^Tq_j\:] =  Q_{:,1:j-1}^TQ_{:,j-1:j}$ }
     \State{$R_{j-1,j-1} = \|q_{j-1}^{j-2}\|_2$, $q_{j-1} = q_{j-1}/R_{j-1,j-1}$ }
     \State{$R_{1:j-1,j} = Q_{j-1}^Tq_j/R_{j-1,j-1}$}
     \State{$T_{1:j-2,j-1} = T(1:j-2,j-1)/R_{j-1,j-1}$ }
     \State{$T_{1:j-2,j-1} = - T_{1:j-2,1:j-2} \times T_{1:j-2,j-1}$ }
     \State{$R_{1:j-1,j} = T_{1:j-1,1:j-1}^T\:R_{1:j-1,j}$ }
     \State{$q_j = q_j - Q_{:,1:j-1} \: R_{1:j-1,j}$ }
\end{algorithmic}
\caption{\footnotesize \label{fig:LGSA} Level-2 Modified Gram-Schmidt Algorithm with $R$ Normalization Lag}
\end{algorithm}

The single global reduction in Algorithm \ref{fig:LGSA} consists of the mass inner product 
in Step 2 and the norm in Step 3. The derivation above can be extended to the iterated
classical Gram-Schmidt algorithm with two full reorthogonalization steps. The columns
of $R$ are updated in a single operation requiring one synchronization and then 
the projection step takes the following form with a second global reduction.
\[
q_j = q_j - Q_{:,1:j-1} \: Q_{:,1:j-1}^T\: (\: q_j - Q_{:,1:j-1} \: R_{1:j-1,j} \:)
\]
\begin{algorithm}
\footnotesize
%\scriptsize
\centering
\begin{algorithmic}[1]
%\State{$\left[ \begin{array}{ccc} Q_{1:m,1:j}, & R_{1:j,1:j}, & T_{1:j,1:j} \end{array} \right] 
\State{${\rm cgs2\_lvl2}
\left[ \begin{array}{ccc} Q_{1:m,1:j}, &  R_{1:j,1:j}, & T_{1:j,1:j} \end{array} \right]$}
     \State{$ [\: Q_{:,1:j-1}^Tq_{j-1},\: Q_{j-1}^Tq_j\:] =  Q_{:,1:j-1}^TQ_{:,j-1:j}$ }
     \State{$R_{j-1,j-1} = \|q_{j-1}^{j-2}\|_2$, $q_{j-1} = q_{j-1}/R_{j-1,j-1}$ }
     \State{$R_{1:j-1,j} = Q_{j-1}^Tq_j/R_{j-1,j-1}$} 
     \State{$T_{1:j-2,j-1} = T(1:j-2,j-1)/R_{j-1,j-1}$ }
     \State{$T_{1:j-2,j-1} = - T_{1:j-2,j-1}$ }
     \State{$R_{1:j-1,j} = T_{1:j-2,j-1}\:R_{1:j-1,j}$ }
     \State{$q_j = q_j - Q_{:,1:j-1} \: Q_{:,1:j-1}^T\: (\: q_j - Q_{:,1:j-1} \: R_{1:j-1,j} \:)$ }
\end{algorithmic}
\caption{\footnotesize \label{fig:LCGSA} Level-2 Iterated Classical Gram-Schmidt Algorithm with $R$ Normalization Lag}
\end{algorithm}
The resulting algorithm only requires two global communication steps.  The
first synchronization step consists of the mass inner product in Step 2 and
normalization in Step 3.  The second synchronization occurs in Step 8 with a
second mass inner product and mass AXPY.  Algorithm \ref{fig:LCGSA} can achieve
orthogonality of the columns of the matrix $Q$ to machine precision.

\section{One Sychronization MGS-GMRES Algorithm}

Given a large sparse linear system of equations $Ax=b$, with residual vector $r
= b - Ax$, the Generalized Minimal Residual (GMRES) algorithm of Saad and
Schultz \cite{Saad86} employs the Arnoldi Gram-Schmidt algorithm with $v_1 =
r_0/\beta$, $\beta= \|r_0\|_2$, in order to orthogonalize the vectors spanning
the Krylov subspace
\[
{\cal K}_m = \{\:v_1,\:A v_1,\:A^2v_1, \ldots A^{m-1}v_1\:\}
\]
The algorithm produces an orthogonal matrix $V_m$ and the 
matrix decomposition $V_m^TAV_m = H_m$
\[
AV_m = V_mH_m + h_{m+1,m} v_{m+1} e_{m+1}^T = V_{m+1} \bar{H}_m 
\]
Indeed, it was observed by Paige et.~al.~\cite{Paige06} that
the algorithm can be viewed as the Gram-Schmidt $QR$ factorization of 
a matrix formed by adding a column to $V_m$ each iteration
\[
\left[ \begin{array}{cc} r_0, & AV_m \end{array} \right] = 
V_{m+1}
\left[ \begin{array}{cc} \|r_0\|e_1, & H_{m+1,m} \end{array} \right]
\]

The GMRES algorithm with modified Gram-Schmidt orthogonalization was derived
by Saad and Schultz \cite{Saad86} and presented as Algorithm \ref{fig:GMRES} below.
\begin{algorithm}
\footnotesize
\centering
\begin{algorithmic}[1]
  \State{$r_1 = b - Ax_1$, $v_1 = r_1 / \|r_1\|_2$.}
  \For {$i=1, 2, \ldots$}
     \State{$z = Av_i$ }
     \For {$j=1,\ldots,i$ }
     \State{$h_{j,i} = (\:z,\:v_j\:)$}
     \State{$w = z -  h_{j,i}v_j$}
     \EndFor
     \State{$h_{i+1,i} = \|w\|_2$ }
     \State{$v_{i+1} = w/h_{i+1,i}$ }
     \State{Apply Givens rotations to $H_{:,i}$ }
  \EndFor
  \State{$y_m = {\rm argmin} \|(\: H_{m+1,m}y_m - \|r_1\|_2e_1 \:) \|_2$ }
  \State{$x = x_1 + V_my_m$ }
\end{algorithmic}
\caption{\footnotesize \label{fig:GMRES} MGS-GMRES}
\end{algorithm}
The standard level-1 BLAS modified Gram-Schmidt (MGS) Algorithm \ref{fig:GMRES}
computes the orthogonalization in Steps 4 to 7.
%
%\[
%\left[ \begin{array}{cc} V_{1:m,i+1}, & H_{1:i+1,i} \end{array} \right] = 
%= {\rm mgs\_lvl1} 
%\left[ \begin{array}{cc} V_{1:m,1:j-1}, &  AV_{1:m,j} \end{array} \right]
%\]
% 
The normalization Steps 8 and 9 are performed afterwards, looking forward
one iteration to compute $h_{i+1,i}$ and $v_{i+1}$
Ghysels et.~al.~\cite{Ghysels13} proposed an alternative form of the
normalization $\|w\|_2$ to avoid a global reduction 
\[
h_{i+1,i} = \left\{\:
\|z\|_2^2 - \|H_{1:i,i}\|_2^2
\:
\right\}^{1/2}
\]
This computation is unstable and suffers from numerical cancellation because it
does not account for the loss of orthogonality in $V$ and approximates the 
inner product
\[
h_{i+1,i}^2 = (\: z - V_{:,i}\:H_{1:i,i} \:)^T(\: z - V_{:,i}\:H_{1:i,i} \:)
\]
Both Ghysels et.~al.~\cite{Ghysels13} and Yamazaki et.~al.~\cite{Yamazaki06}
describe pipelined algorithms based on CGS-1 that incorporate an iteration lag
in the normalization step to avoid the unstable normalization above. However,
these algorithms are somewhat unstable because the loss of orthogonality for
CGS-1 is ${\cal O}(\eps)\kappa(A)^2$, see \cite{Giraud06}. The authors employ a
change of basis $QR$ factorization, $V_k = Z_kG_k$ and shifts in order to
mitigate these instabilities. However, the pipeline depths still must remain
relatively small.

The MGS-GMRES algorithm is backward stable \cite{Paige06} and orthogonality is
maintained to ${\cal O}(\eps)\kappa(A)$, where $\kappa(A)$ is the condition
number of $A$ defined as the ratio of the largest to smallest singular values
\cite{Giraud05}.  By employing the Level-2 BLAS compact $WY$ MGS algorithm with
a backward lag, the normalization for the $i$--th iteration is computed to full
accuracy with the orthogonalization for iteration $i+1$. In effect, the lag
creates a pipeline where orthogonalization of the Krylov vectors $v_1$, and
$v_2$ is performed before the diagonal of the $R$ matrix in the $QR$
factorization, corresponding to $h_{2,1}$, is computed for iteration 2.  The
lagged compact $WY$ algorithm is implemented below as Algorithm
\ref{fig:onesynch} for the $i$-th iteration of MGS-GMRES, where $V_{1:m,i+2} = AV_{1:m,i+1}$.
%
%\[
%\left[ \begin{array}{ccc} V_{1:m,1:i+2}, & R_{1:i+2,1:i+2}, & T_{1:i+2,1:i+2} \end{array} \right] 
%= {\rm mgs\_lvl2} 
%\left[ \begin{array}{ccc} V_{1:m,1:i+2}, &  R_{1:i+2,1:i+2}, & T_{1:i+2,1:i+2} \end{array} \right]
%\]
%
\begin{algorithm}
\footnotesize
\centering
\begin{algorithmic}[1]
   \State{$r_1 = b - Ax_1$, $v_1 = r_1 / \|r_1\|_2$.}
   %\State{$\left[\: V_{1:m,1:1},\: R_{1:1,1:1},\: T_{1:1,1:1}\: \right] = 
   \State{${\rm mgs\_lvl2} 
           \left[\: V_{1:m,1:1},\: R_{1:1,1:1},\: T_{1:1,1:1}\: \right]$ }
   \State{$v_{2} = Av_{1}$ } 
   %\State{$\left[\: V_{1:m,1:2},\: R_{1:2,1:2},\: T_{1:2,1:2}\: \right] = 
   \State{${\rm mgs\_lvl2} 
           \left[\: V_{1:m,1:2},\: R_{1:2,1:2},\: T_{1:2,1:2}\: \right]$ }
   \State{$R_{2,2} = \|v_{2}\|_2$, $v_{2} = v_{2} / R_{2,2}$ } 
   \State{$H_{1:2,1} = R_{1:2,2}$}
   \State{Apply Givens rotations to $H_{:,1}$ }
  \For {$i=1, 2, \ldots$}
     \State{$v_{i+2} = Av_{i+1}$ } 
     %\State{$\left[\: V_{1:m,1:i+2},\: R_{1:i+2,1:i+2},\: T_{1:i+2,1:i+2}\: \right] = 
     \State{${\rm mgs\_lvl2} 
             \left[\: V_{1:m,1:i+2},\: R_{1:i+2,1:i+2},\: T_{1:i+2,1:i+2}\: \right]$ }
     \State{$R_{i+2,i+2} = \|v_{i+2}\|_2$, $v_{i+2} = v_{i+2} / R_{i+2,i+2}$ } 
     \State{$H_{1:i+2,i+1} = R_{1:i+2,i+2}$}
     \State{Apply Givens rotations to $H_{:,i+1}$ }
  \EndFor
  \State{$y_m = {\rm argmin} \|(\: H_{m+1,m}y_m - \|r_1\|_2e_1 \:) \|_2$ }
  \State{$x = x_1 + V_my_m$ }            
\end{algorithmic}
\caption{\footnotesize \label{fig:onesynch} One Synchronization MGS-GMRES with Lagged Normalization}
\end{algorithm}
The normalization for iteration $i+2$ is computed in Steps 11 and 12.
%
%\[
%\begin{array}{ccc} 
%R_{i+2,i+2} = \|v_{i+2}\|_2, & v_{i+2} = v_{i+2} / R_{i+2,i+2}, & H_{1:i+2,i+1} = R_{1:i+2,i+2}
%\end{array} 
%\]
%
Therefore, the column $H_{1:i+2,i+1} = R_{1:i+2,i+2}$ of the Hessenberg matrix
is pre-computed for iteration $i+1$ and the matrix $H_{:,i}$ is employed to
solve the least squares problem at iteration $i$. Because the loss of 
orthogonality for MGS-GMRES is ${\cal O}(\eps)\kappa(A)$, Algorithm \ref{fig:onesynch}
is not only more stable but can be pipelined to an arbitrary iteration depth $L$.

%%In order to illustrate the creation of a pipeline, consider the example below with 
%%a pipeline of depth $L=2$.
%
\begin{algorithm}
\footnotesize
\centering
\begin{algorithmic}[1]
   \State{$r_1 = b - Ax_1$, $v_1 = r_1 / \|r_1\|_2$.}
   \State{${\rm mgs\_lvl2} \left[\: V_{1:m,1:1},\: R_{1:1,1:1},\: T_{1:1,1:1}\: \right]$ }
   \State{$v_{2} = Av_{1}$ } 
   \State{${\rm mgs\_lvl2} \left[\: V_{1:m,1:2},\: R_{1:2,1:2},\: T_{1:2,1:2}\: \right]$ }
   \State{$R_{2,2} = \|v_{2}\|_2$, $v_{2} = v_{2} / R_{2,2}$ } 
   \State{$H_{1:2,1} = R_{1:2,2}$}
   \State{Apply Givens rotations to $H_{:,1}$ }
  \For {$i=1, 3, \ldots$}
     \State{$v_{i+2} = Av_{i+1}$ } 
     \State{${\rm mgs\_lvl2} \left[\: V_{1:m,1:i+2},\: R_{1:i+2,1:i+2},\: T_{1:i+2,1:i+2}\: \right]$ }
     \State{$R_{i+2,i+2} = \|v_{i+2}\|_2$, $v_{i+2} = v_{i+2} / R_{i+2,i+2}$ } 
     \State{$H_{1:i+2,i+1} = R_{1:i+2,i+2}$}
     \State{Second orthogonalization}
     \State{$v_{i+3} = Av_{i+2}$ } 
     \State{${\rm mgs\_lvl2} \left[\: V_{1:m,1:i+3},\: R_{1:i+3,1:i+3},\: T_{1:i+3,1:i+3}\: \right]$ }
     \State{$R_{i+3,i+3} = \|v_{i+3}\|_2$, $v_{i+3} = v_{i+3} / R_{i+3,i+3}$ } 
     \State{$H_{1:i+3,i+2} = R_{1:i+3,i+3}$}
     \State{Apply Givens rotations to $H_{:,i+1}$ }
     \State{Apply Givens rotations to $H_{:,i+2}$ }
  \EndFor
\end{algorithmic}
\caption{\footnotesize \label{fig:pipeline} $L=2$ MGS-GMRES Pipeline}
\end{algorithm}

\section{Implementation Details}

The standard Level-1 BLAS MGS-GMRES implementation with the modified Gram-Schmidt process
requires $i+1$ separate inner-products (global summations) in iteration $i$.
Therefore, the amount of global communication grows quadratically.  Several
authors have proposed how to reduce the communication overhead associated with
MGS-GMRES.  The iterated classical Gram-Schmidt algorithm
with two reorthogonalization steps requires at least two global reductions. The
Trilinos-Belos solver stack employs an iterated classical Gram-Schmidt
ICGS-GMRES \cite{Belos}.  There have been several
approaches proposed in the literature to reduce the communication overhead
associated with ICGS-GMRES.
%the classical Gram-Schmidt Arnoldi algorithm in 
Hernandez et.~al.~\cite{Hernandez} implement the re-orthogonalization with
communication pipelining.  The approach developed here requires less
computation and communication than previous algorithms.  The one-sync Level-2 BLAS
MGS-GMRES algorithm allows us to perform the compact $WY$ modified Gram-Schmidt algorithm,
by combining a mass inner-product and normalization into one global reduction
{\tt MPI\_AllReduce} per GMRES iteration.  The mass inner product for computing
$v = X^Ty$ is written as 
\[
v = y \cdot X = y\cdot x_1 + y\cdot x_2 + \cdots + y\cdot x_m
\]
The algorithm for $w = z - V\:H_{1:i,i}$ is implemented with the vector mass AXPY (MAXPY).
\[
y = y + X\begin{bmatrix}\alpha_1 &&&&\\
&\alpha_2&&&\\
&&&&\\
&&&&\\
&&&&\alpha_m\end{bmatrix} = y + \alpha_1x_1 + \ldots + \alpha_m x_m
\]
The PetSC library employs the vector MAXPY for a pipelined $p^1$--GMRES iteration 
with reorthogonalization \cite{Ghysels13}.

For optimal performance of these two algorithms on current generation
multi-core processors, cache-blocking combined with loop unrolling can be
exploited. Once the Krylov sub-space dimension is sufficiently large, the outer
loop of the mass inner-product should be unrolled into batches of two, four, or
more summations.  The overall effect is to expose more work to multiple
floating point units and pre-fetch $y$ together with multiple columns of $X$
into the cache.  The outer loop can also be threaded given the typically large
vector lengths.  In an MPI parallel implementation, the columns of $X$ are
partitioned across the MPI ranks. Unrolling should have the beneficial effect
of increasing the memory band-width utilization of the processors. The vector
MAXPY should also be unrolled with the elements $\alpha_i$ pre-fetched into the
cache. The Hessenberg matrix $H$ is small and dense. Thus, it is laid out in
memory by columns for cache access and fast application of the orthogonal
rotations in the least-squares problem.  The remaining operations inside
the one-sync MGS-GMRES are to compute norms, which are BLAS-1 operations, (sparse)
matrix-vector products and the application of a preconditioner. BLAS-1
operations are usually highly optimized for the particular architecture. The
sparse matrix-vector multiply is parallel and can be easily partitioned by rows
when using a compressed-sparse row (CSR) storage format.

Both the mass inner product and vector MAXPY are well-suited to a
single-instruction multiple-data (SIMD) model of parallel execution. These
computations can be massively threaded on a GPU. For GPU execution, the mass
inner-product can be off-loaded to the GPU but requires multiple CUDA kernel
launches to perform synchronizations between reduction steps. Synchronization is
possible within streaming multi-processors (SM) executing thread blocks.
However, a CUDA kernel launch is needed to synchronize across the SM.
Unrolling of the inner-products into batches is also possible and permits
higher sustained memory bandwidth between the GPU global memory and streaming
multi-processors (SM). This is also the case for the vector MAXPY.  A CUDA
implementation of the vector MAXPY is provided by PetSC. A recent paper by
Romero et.~al.~\cite{Romero14} also describes a vector MAXPY for the GPU.  
The GPU implementation of the CGS-1 based CA-GMRES algorithm is described in \cite{Yamazaki14}.
An important optimization to save storage and promote data re-use is to
fuse the matvec+precon step and not store the preconditioned vectors \cite{SwiryThesis}.

\section{Numerical Results}

The primary motivation for proposing the low synchronization GMRES algorithms
was to improve the performance and strong-scaling characteristics of DOE
physics based simulations at exascale and, in particular, the Nalu CFD solver,
Domino et.~al.~\cite{Domino17}. The Nalu model solves the incompressible
Navier-Stokes equations on unstructured finite-element meshes. The momentum and
pressure equations are time-split and solved sequentially. Both the momentum
and pressure continuity equation are solved with preconditioned GMRES
iterations.  The pressure preconditioner employs either smoothed-aggregation or
classical Ruge-St\"uben algebraic multigrid (C-AMG). The latter is provided by
the Hypre-BoomerAMG library of solvers \cite{Henson00}.  The Krylov solvers,
which typically involve a substantial amount of communication, are an expensive
part of the simulation, requiring over 50\% of the simulation time.

In order to evaluate the performance of the different GMRES solvers, a separate
Hypre linear solver application has been implemented.  In all our performance
studies, matrices were obtained from Nalu wind turbine simulations and tested
in our Hypre solver application.  The performance studies involving GPU
acceleration were completed on a single node of the Summit-Dev supercomputer at
ORNL.  Summit-Dev is built with IBM S822LC nodes which consist of two 10-core
IBM Power-8 processors with 256 GB of DD4 memory and four NVIDIA Tesla P100
GPUs connected by NVlink running at 80 GB/sec. The GPUs share 16 GB HBM2
memory. The interconnect is based on two rails of Mellanox EDR Infiniband.

Gram-Schmidt orthogonalization, as emphasized earlier, is the most expensive part
of the GMRES algorithm. Hence, the algorithms developed in this paper are applied.
These are based upon one and two synchronization CGS and MGS Gram-Schmidt algorithms.
For each example, we specify which orthogonalization method was employed.

\subsection{Numerical Stability and Accuracy}

\begin{figure}
\includegraphics[width=0.45\textwidth]{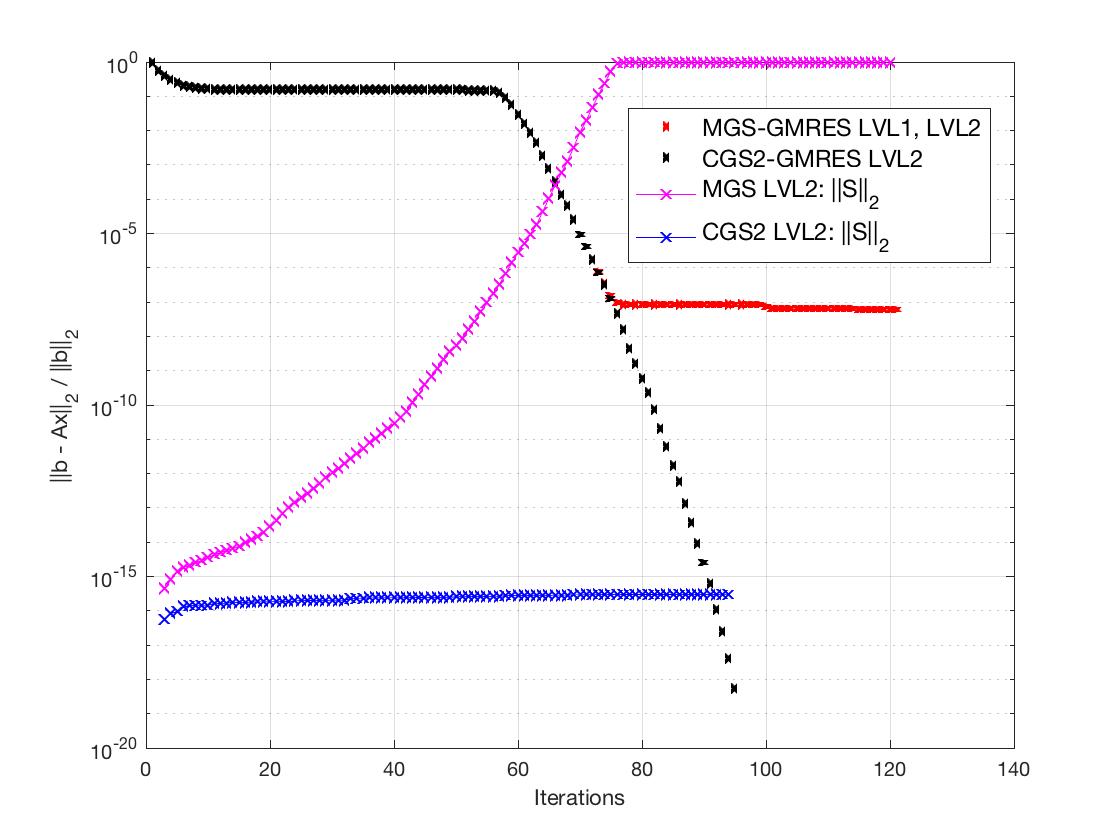}
\caption{Simonici matrix. Comparison of limiting accuracy for MGS-GMRES and CGS2-GMRES.
Loss of orthogonality measure $\|S\|_2$ \label{fig:m1}}
\end{figure}

Linear systems appearing in the literature are employed to evaluate the
numerical stability and limiting accuracy for several of the GMRES algorithms
described here. All tests were performed in Matlab and the loss of
orthogonality metric $\|S\|_2$ that is plotted for the Krylov vectors was
derived by Paige \cite{Paige18}. The first problem was proposed by Simonici and
Szyld \cite{Simonici03}. A diagonal matrix of dimension 100 is specified as $A
= {\rm diag}(1e-8,2, \ldots, 100)$ with random right-hand side $b = {\rm
randn}(100,1)$, that is normalized so that $\|b\|_2 =1$.  The matrix condition
number is $\kappa(A) = 1\times 10^{10}$ and is controlled by the small first
diagonal element. A preconditioner is not employed in this case.
The problem allows us to specify the limiting accuracy of the
GMRES iteration. In particular, Figure \ref{fig:m1} plots the relative residual
versus the number of iterations for the standard level-1 BLAS MGS-GMRES and
level-2 one-synch lagged MGS-GMRES algorithm.  These exhibit identical
convergence behaviour with the loss of orthogonality metric $\|S\|_2$
increasing to one at 80 iterations. The convergence stalls at this point with
the relative residual reaching $1\times 10^{-7}$.  In the case of the one or
two-synch classical Gram-Schmidt (CGS-2) based on Algorithm 5, the loss
of orthogonality metric $\|S\|_2$ remains close to machine precision and the
convergence curve continues to decrease without stalling to $1\times 10^{-18}$.

\begin{figure}
\includegraphics[width=0.45\textwidth]{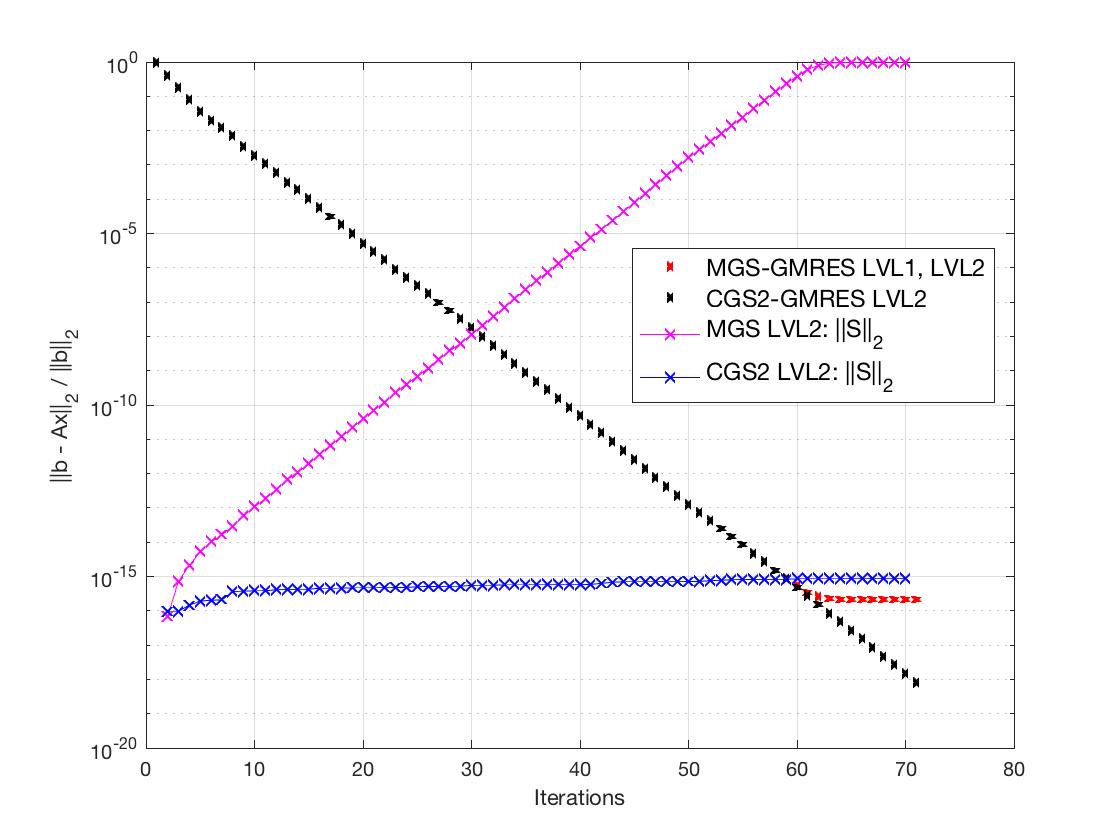}
\caption{Thermal1  matrix. Comparison of limiting accuracy for MGS-GMRES and CGS2-GMRES.
Loss of orthogonality measure $\|S\|_2$ \label{fig:m2}}
\end{figure}

A second problem is taken from the Florida (now Texas A\&M) sparse matrix
collection. The particular system is `thermal1' with dimension $82654$
and number of non-zeros, $nnz(A) = 574458$. The right-hand side is
obtained from $b=Ax$, where $x = [\:1,\:1,\:\ldots,\:1\:]^T$. For this
problem a Ruge-St\"uben \cite{Ruge87} classical AMG preconditioner is applied.
Once again the standard level-1 MGS-GMRES and level-2 one-synch GMRES
algorithms are compared with the CGS-2 GMRES algorithm(s). The relative
residuals and loss of orthogonality metrics are plotted in Figure \ref{fig:m2}.
Observe that when $\|S\|_2 = 1$ that the convergence stalls for the MGS-GMRES
algorithms, whereas the CGS-2 variants maintain orthogonality to near the
level of machine precision.

\subsection{GPU C/CUDA Performance}

The results presented below were obtained on Summit-Dev computational cluster
(OLCF) using NVCC V9.0.67. All variables are allocated using device memory.
The standard cuSPARSE and cuBLAS functions are used to perform BLAS operations
(SpMV, AXPY, SCAL, DOT). It could be assumed that MDOT and MAXPY can be
expressed as dense matrix-matrix product (MDOT) and a combination of dense
matrix-matrix products (MAXPY), and implemented using standard cuBLAS
routines. However, it was found that this approach was less efficient than 
our custom implementation. We speculate that the main reason for sub-optimal
performance was the need to transpose $V$ in MDOT and to perform two CUDA
kernel calls instead of one for MAXPY. Hence, MAXPY and MDOT are implemented
as fused operations in CUDA and optimized for maximum performance. 

The problem employed for GPU performance analysis was extracted from the Nalu
Vestas V27 wind turbine simulation described in Thomas
et.~al.~\cite{Thomas18}. The discretization scheme is based on the
Control-Volume Finite-Element method (CVFEM). The matrix has $675905$ rows and
$675905$ columns. The average number of non-zeros per row is $11$. The matrix
is stored in the standard compressed-sparse row (CSR) format, which is
supported by the cuSPARSE library.  A preconditioner is not applied in this test.

\begin{figure}
\includegraphics[width=0.45\textwidth]{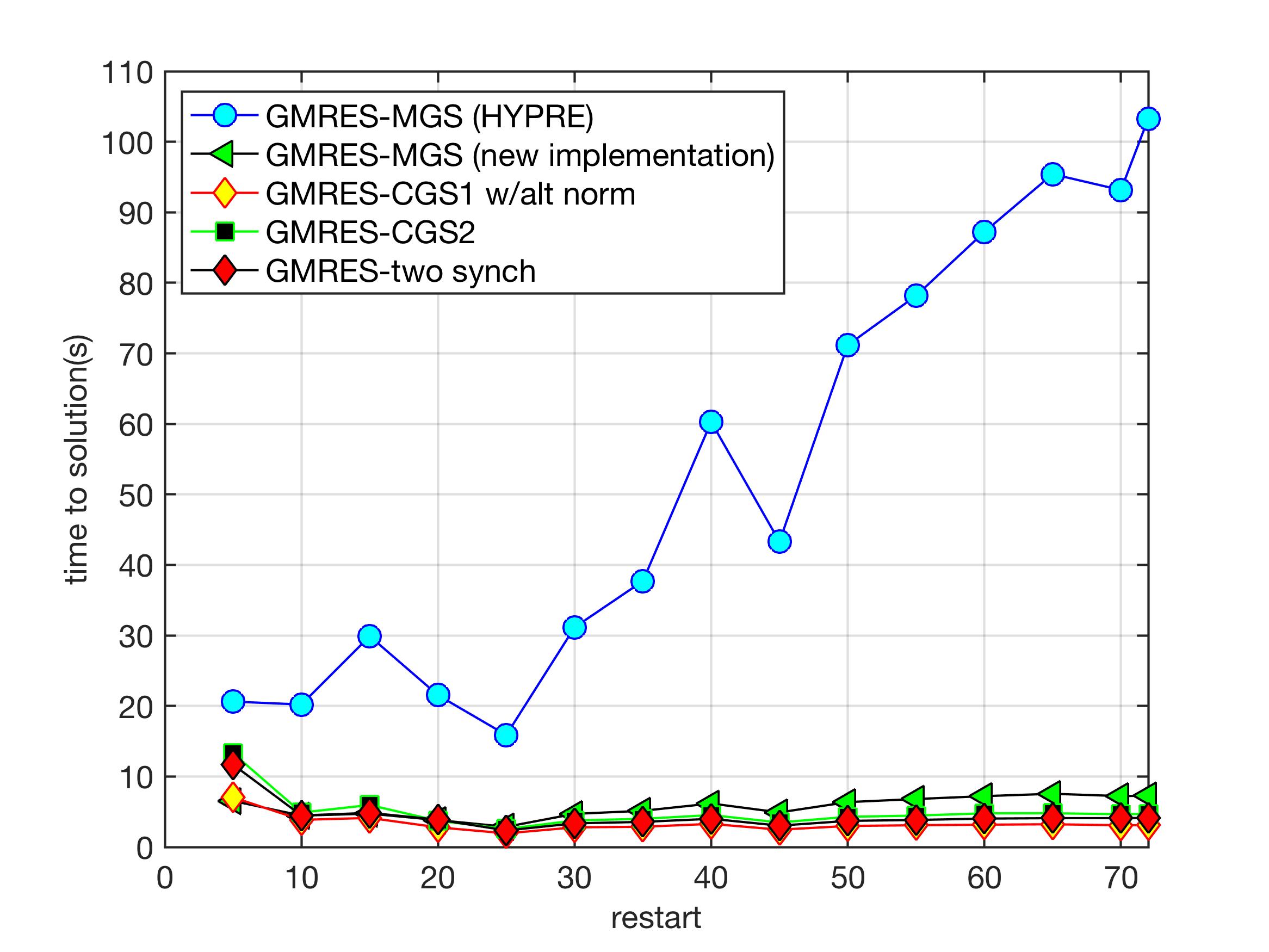}
\caption{Performance of the Hypre GPU implementation of Level-1 BLAS MGS-GMRES(m) versus
GMRES(m) with alternative orthogonalization strategies. Tolerance $1\times 10^{-6}$.
Linear system from Nalu wind simulation\label{fig:fig0}}
\end{figure}

Five approaches are compared for solving the linear system: (1) Level-1 BLAS
MGS-GMRES(m) compiled using CUDA, (2) MGS-GMRES(m) written by the authors for
Hypre using only device memory and reducing unnecessary data copies and
convergence checks, (3) GMRES(m) with CGS-1 and Ghysels normalization, (4)
CGS-2 GMRES (Algorithm 1), and (5) two-synch GMRES(m) as in Algorithm 3.
Restarts are set as $m = 5, 10, \ldots, 65, 70, 72$ with relative tolerance
$1\times 10^{-6}$.  The only difference between approaches (2)--(5)  is how
they orthogonalize Krylov vectors.  Figure~\ref{fig:fig0} displays the improvement
in run-time by using alternative orthogonalization approaches. The improvement
is much greater, as expected, for large restart values.  Figure~\ref{fig:fig2}
indicates the associated speedups.

\begin{figure}
\includegraphics[width=0.45\textwidth]{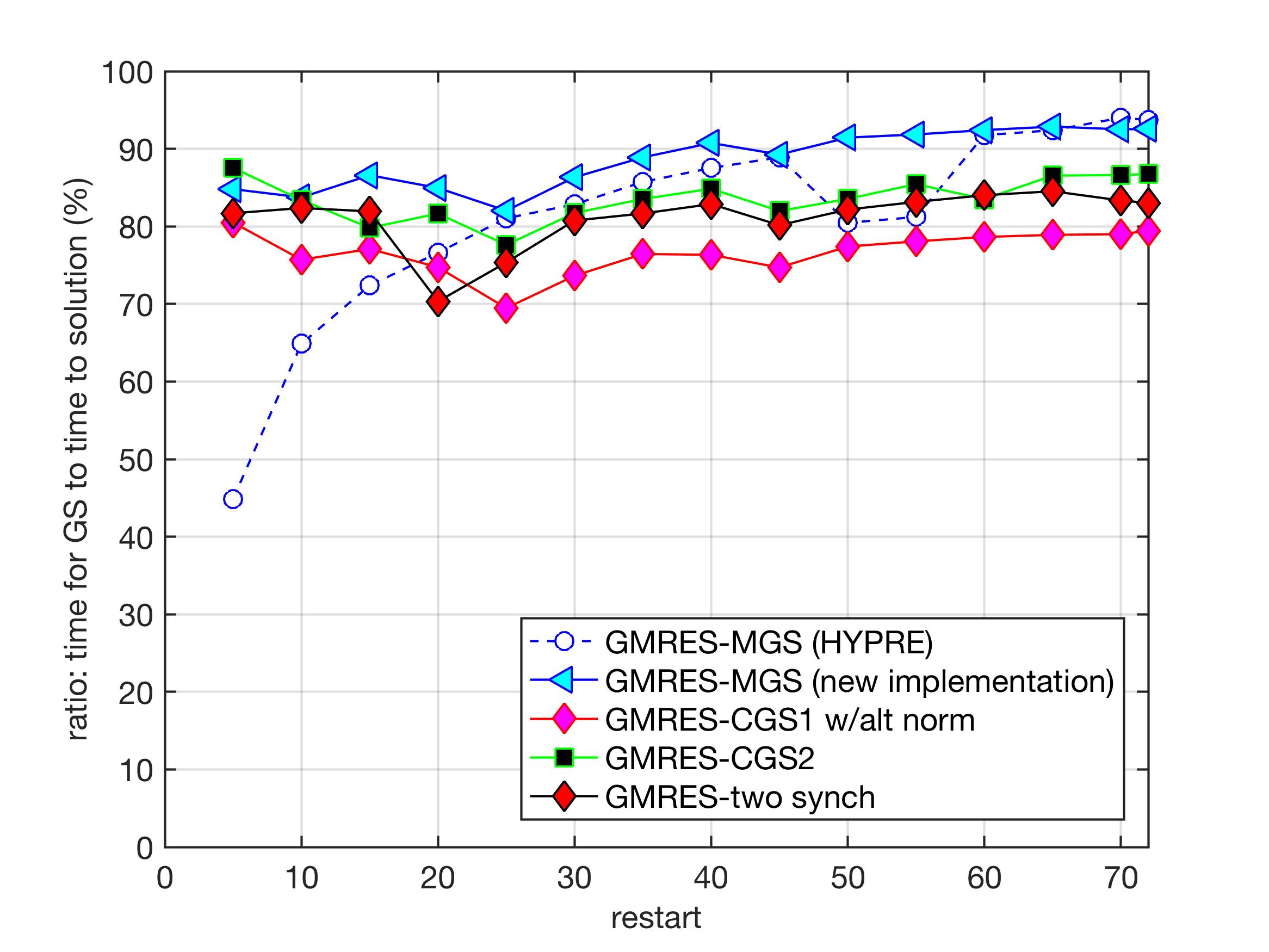}
\caption{Performance of the Hypre GPU implementation of Level-1 BLAS MGS-GMRES(m) with MGS and
GMRES(m) with alternative orthogonalization strategies. The ratio of the time
spent on Gram-Schmidt orthogonalization versus time to solution.
\label{fig:fig1}}
\end{figure}

\begin{figure}
\includegraphics[width=0.45\textwidth]{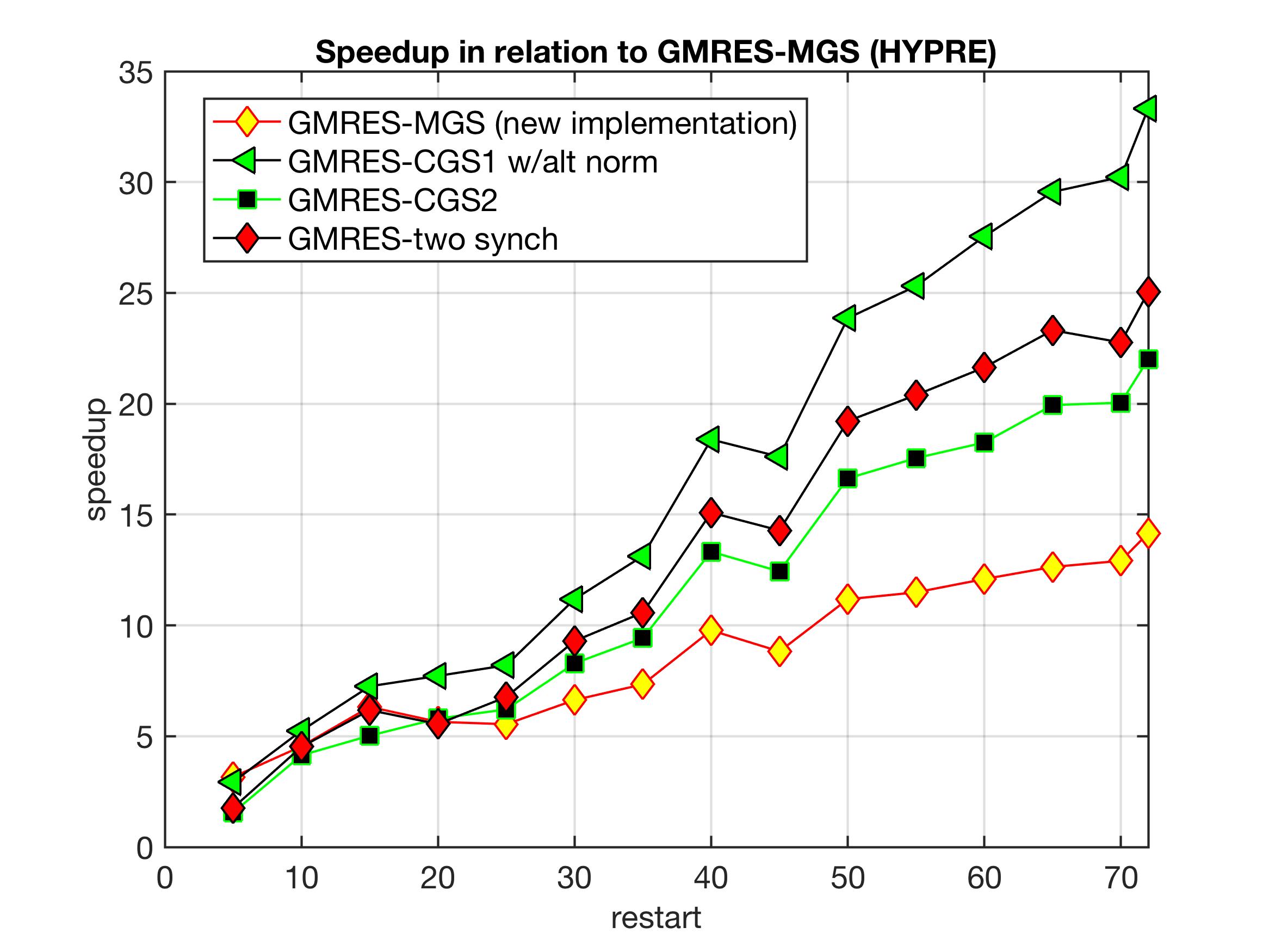}
\caption{ The speedup resulting from replacing Hypre Level-1 BLAS modified Gram-Schmidt with
alternate low-latency GMRES variants.  \label{fig:fig2}}
\end{figure}

Figure~\ref{fig:fig1} is a plot of the ratio of wall-clock run-time taken by
Gram-Schmidt to the total run-time of the sparse linear system solver. One can
easily observe that the Level-1 BLAS modified Gram-Schmidt accounts for 90\% of
the total run-time. Using GMRES with different orthogonalization strategies
lowers this time down to 40\% from 70\%. The remainder of the execution time is
taken by various BLAS routines and matvecs. Hence, an even larger improvement
is expected in run-time when using a preconditioner.
Figure~\ref{fig:fig2} displays the speedup resulting from using low-latency
orthogonalization methods. The implementation of the Hypre MGS-GMRES
GPU implementation is taken as the baseline reference timing.

\begin{figure}
\includegraphics[width=0.45\textwidth]{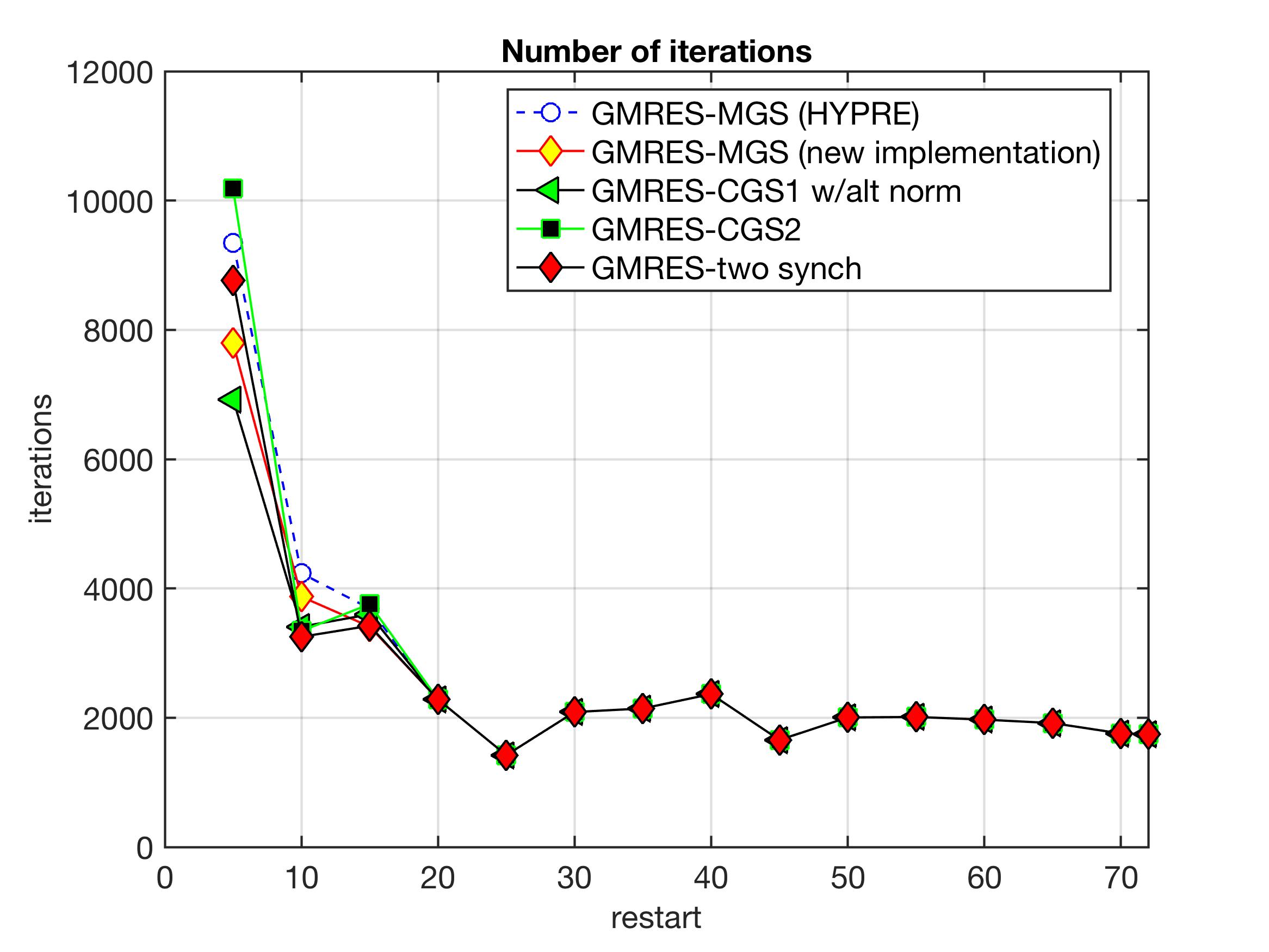}
\caption{Number of iterations for Hypre Level-1 BLAS MGS-GMRES(m) versus alternate
low latency GMRES(m). \label{fig:fig3}}
\end{figure}

Figure~\ref{fig:fig3} displays the number of GMRES iterations resulting from using 
low-latency orthogonalization methods. It may be observed that for short restart 
cycles the number of iterations varies for the different methods. This is expected 
because of the CGS-1 stability versus CGS-2 and MGS.

The `Parabolic FEM' from the Florida sparse matrix collection is used in a
second performance test. Total run-times are reported in Table
\ref{table:solve}.  The matrix dimension is $n=525825$ with $nnz(A) = 3674625$
non-zeros.  The Hypre-BoomerAMG $V$-cycle with 3 levels is the preconditioner
with $L_1$ scaled Jacobi smoother on the GPU. The total solver (Solve) time is
reported along with the Gram-Schmidt (Orth) times and number of GMRES
iterations (Iters) required to reach the specified relative tolerance. The
lowest run-times are achieved by the GMRES two-synch algorithm.  
These are 60\% faster than the Hypre Level-1 BLAS MGS-GMRES (for $1\times 10^{-13}$).

\begin{table}
\renewcommand{\arraystretch}{1.5}
\begin{tabular}{|l|l|l|l|l|} \hline
  Method & tolerance & Iters & Orth (s) & Solve (s) \\
  \hline
  GMRES-HYPRE         & 1e-13 & 78  & 0.95  & 3.06 \\
  GMRES-MGS (new)     & 1e-13 & 78 & 0.25  & 2.03 \\
  GMRES-CGS1 alt norm & 1e-13 & N/A & N/A   & N/A  \\
  GMRES-CGS2          & 1e-13 & 78 & 0.14  & 2.00 \\
  GMRES-two synch     & 1e-13 & 78 & 0.12  & 1.92 \\ \hline
  GMRES-HYPRE         & 1e-15 & N/A & N/A   &  N/A  \\
  GMRES-MGS (new)     & 1e-15 & 93 & 0.27 & 2.34 \\
  GMRES-CGS1 alt norm & 1e-15 & N/A & N/A   & N/A   \\
  GMRES-CGS2          & 1e-15 & 92 &  0.15 &  2.14 \\
  GMRES-two synch     & 1e-15 & 92 & 0.13  & 2.26  \\  \hline
\end{tabular}
\caption{Parabolic FEM matrix, restart: $72$, N/A - did not converge. \label{table:solve} }
\end{table}

\begin{figure}
\centering
\includegraphics[width=0.45\textwidth]{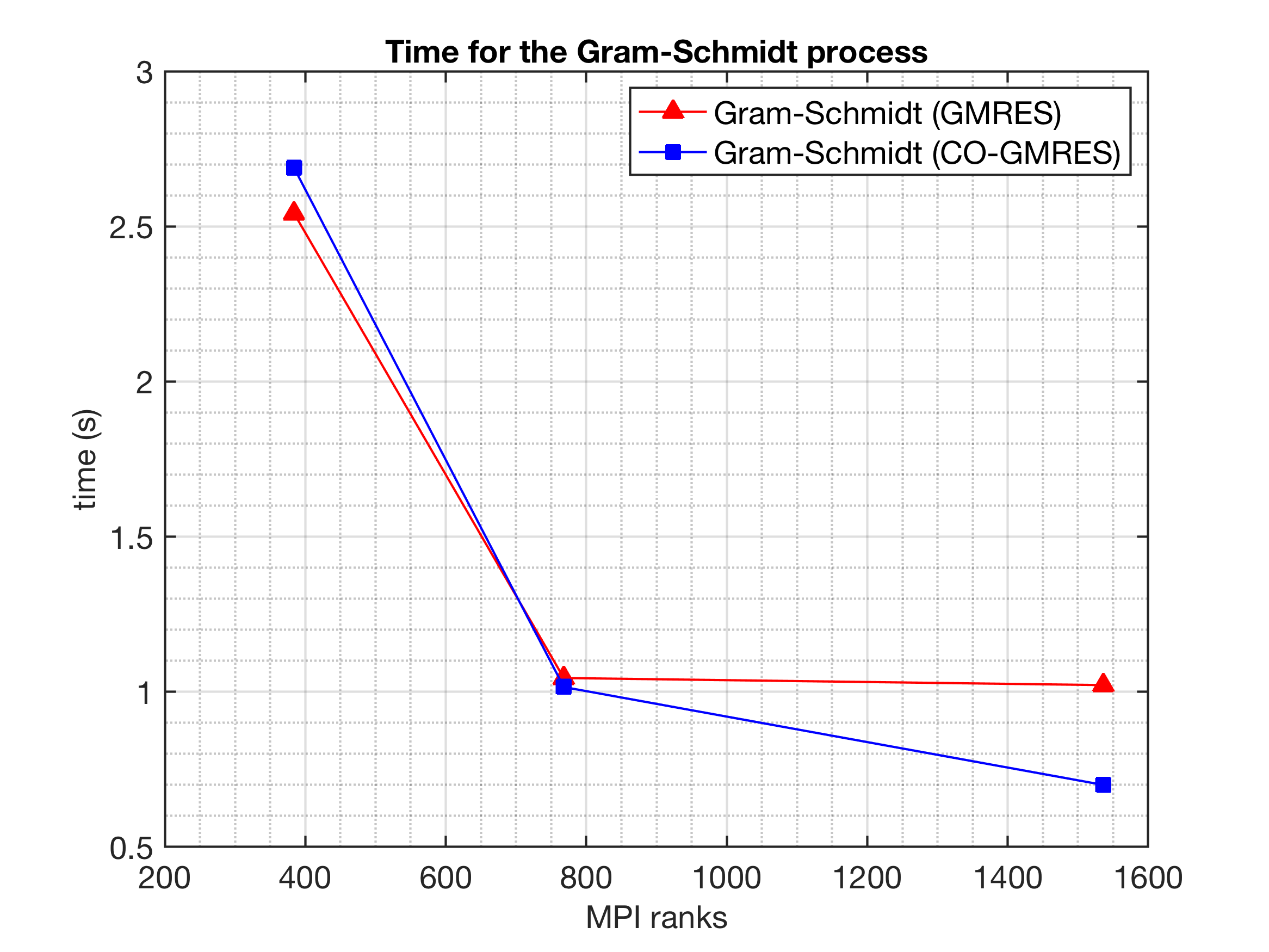}\\
\caption{Gram-Schmidt kernels time on Peregrine at NREL.
Time for a single linear system solve. Matrix obtained from
Vestas V27 66M element mesh with 29M DOF.}
\label{fig:GS}
\end{figure}

Initial testing of the Hypre linear solver on the Peregrine supercomputer at
NREL has demonstrated that the low-latency algorithms result in a reduced
communication overhead for the MPI global reductions.  Peregrine nodes contain Xeon
E5-2670 v3 Haswell processors, $2\times 12$ core sockets,  64 GB DDR4 memory,
inter-connected with an Infiniband network.  The run-time spent in the Gram-Schmidt
orthogonalization is plotted in Figure \ref{fig:GS} for the Nalu V27 linear system.
Clearly, the new algorithms reduce the time as the node count increases.

\section{Conclusions and Future goals}

We have presented low-latency classical and modified Gram-Schmidt algorithms
and applied these to the Arnoldi-GMRES algorithm. The number of synchronization
steps depends on the desired level of orthogonality, stability and limiting
accuracy of the the resulting GMRES Krylov iterations. A one-synchronization
MGS-GMRES algorithm based on the compact $WY$ MGS algorithm with a lagged
normalization step results in a backward-stable algorithm \cite{Paige06} with
${\cal O}(\eps)\kappa(A)$ loss of orthogonality.  GPU implementations of these
low-latency algorithms achieve up to a 35 times speed-up over the 
Hypre Level-1 BLAS MGS-GMRES algorithm. %A conservative estimate of GPU
%acceleration that is 50 to 100 times faster than the CPU implementation has
%been achieved.

The execution speed of our GPU low-latency GMRES algorithms implies a new
scaling paradigm for large sparse linear solvers. The relative balance in
run-time that is split between GMRES (including matrix-vector multiplies) and a
preconditioner such as the Hypre-BoomerAMG $V$-cycle may now shift in favor of
a larger number of inexpensive GMRES iterations combined with a light-weight
$V$-cycle in order to achieve a much lower and optimal run-time. For example,
we solved a large $n=500$K linear system in under two seconds with 100 GMRES
iterations, using GPU accelerated matrix-vector multiplies, preconditioned by a
less costly $V$-cycle and $L_1$ Jacobi smoother on the GPU.  The challenge
is to extend this scaling result to a large number of nodes based on a low
number of synchronizations.

An obvious extension of this work is a distributed-memory parallel
implementation and testing using MPI+OpenMP and MPI+CUDA. A first step in this
direction would be a multi-GPU implementation. The theory developed in this
paper and the preliminary results indicate that our approach, while stable,
also has the desired scalability properties both in terms of fine and coarse
grain parallelism.  Further, investigation is needed to establish the
MPI-parallel properties of the proposed algorithms at large scale on 
${\cal O}(100)K$ nodes.

\section{Acknowledgements}

This work was authored in part by the National Renewable Energy Laboratory,
operated by Alliance for Sustainable Energy, LLC, for the U.S. Department of
Energy (DOE) under Contract No. DE-AC36-08GO28308. Funding provided by the
Exascale Computing Project (17-SC-20-SC), a collaborative effort of two U.S.
Department of Energy organizations (Office of Science and the National Nuclear
Security Administration) responsible for the planning and preparation of a
capable exascale ecosystem, including software, applications, hardware,
advanced system engineering, and early testbed platforms, in support of the
nation's exascale computing imperative. The views expressed in the article do
not necessarily represent the views of the DOE or the U.S. Government. The U.S.
Government retains and the publisher, by accepting the article for publication,
acknowledges that the U.S. Government retains a nonexclusive, paid-up,
irrevocable, worldwide license to publish or reproduce the published form of
this work, or allow others to do so, for U.S. Government purposes. 

%\section{Appendix: performance optimization of MAXPY and MDOT on the GPU.}
%
%\subsection{MAXPY}
%While performing MAXPY involving $k$ vectors with dimension $n \times 1$, we
%need to copy all the vectors and the scaling factors to the GPU and retrieve
%the $n \times 1$ solution vector from the GPU. Because speed of the data copy
%is a limiting factor for the GPU performance, we assume that the code cannot
%perform any faster than the data copy. This assumption has been applied to the
%GPU performance modeling in \cite{Swirydowicz17}.  One might argue that MAXPY
%is simply a matrix-matrix product, i.e., it can be written as $y - X \alpha$.

\end{document}